\documentclass[number,citesort,dvips]{arxbj}

\aid{0}
\volume{17}
\issue{1}
\pubyear{2011}
\firstpage{211}
\lastpage{225}
\doi{10.3150/10-BEJ267}

\makeatletter

\newproclaim{definition}{Definition}[section]
\newtheorem{theorem}[definition]{Theorem}
\newremark{example}[definition]{Example}
\newtheorem{lemma}[definition]{Lemma}
\newtheorem{proposition}[definition]{Proposition}

\newcommand{\snorm}[1]{\Vert #1 \Vert}
\newcommand{\hnorm}[1]{\Vert #1 \Vert_{H}}
\newcommand{\inorm}[1]{\Vert #1 \Vert_\infty}
\newcommand{\Lx}[2]{{L}_{#1}(#2)}
\newcommand{\clippt} {{{}^{\smallfrown}} \hspace*{-1.5ex} t\hspace*{0.5ex}}
\newcommand{\clippfo} {{{}^{^{\textstyle \smallfrown}}} \hspace*{-2.1ex} f}
\newcommand{\fTc}{{{}^{^{\textstyle \smallfrown}}} \hspace*{-2.1ex} f_{D,\lambda}}
\newcommand{\fTcn}{{{}^{^{\textstyle \smallfrown}}} \hspace*{-2.1ex} f_{D,\lambda_n}}
\newcommand{\RP}[2]{{\mathcal{R}_{#1,\mathrm{P}}(#2)}}
\newcommand{\RPB}[1]{{\mathcal{R}_{#1,\mathrm{P}}^{*}}}
\newcommand{\RT}[2]{{\mathcal{R}_{#1,\mathrm{D}}(#2)}}
\newcommand{\Rx}[3]{{\mathcal{R}_{#1,#2}(#3)}}
\newcommand{\fpb}{{f_{\tau,\mathrm{P}}^*}}
\newcommand{\CQ}[2]{{\mathcal{C}_{#1,\mathrm{Q}}(#2)}}
\newcommand{\CQB}[1]{{\mathcal{C}_{#1,\mathrm{Q}}^{*}}}
\newcommand{\CPo}[1]{{\mathcal{C}_{#1,\mathrm{P}(\cdot|x)}}}
\newcommand{\CPoB}[1]{{\mathcal{C}_{{#1,\mathrm{P}(\cdot|x)}}^{*}}}
\newcommand{\Cxx}[3]{{\mathcal{C}_{#1,#2}(#3)}}
\newcommand{\FQ}[2]{\mathcal{M}_{#1,\mathrm{Q}}(#2)}
\newcommand{\FPo}[2]{\mathcal{M}_{#1,\mathrm{P}(\cdot|x)}(#2)}
\newcommand{\dmaxexu}[3]{{\delta}_{\max,#3}(#1,#2)}
\newcommand{\dmaxexub}[3]{{\delta}_{\max,#3}(#1,#2)}

\makeatother

\begin{document}
\begin{frontmatter}

\title{Estimating conditional quantiles with the help of the pinball loss}
\runtitle{Estimating conditional quantiles with the help of the pinball loss}

\begin{aug}
\author[a]{\fnms{Ingo} \snm{Steinwart}\corref{}\thanksref{a}\ead[label=e1]{ingo.steinwart@mathematik.uni-stuttgart.de}} \and
\author[b]{\fnms{Andreas} \snm{Christmann}\thanksref{b}\ead[label=e2]{Andreas.Christmann@uni-bayreuth.de}}
\runauthor{I. Steinwart and A. Christmann}
\address[a]{University of Stuttgart,
                      Department of Mathematics,
                      D-70569 Stuttgart, Germany.
                      \\
\printead{e1}}
\address[b]{University of Bayreuth,
Department of Mathematics,
D-95440 Bayreuth.\\
\printead{e2}}
\end{aug}

\received{\smonth{11} \syear{2008}}
\revised{\smonth{1} \syear{2010}}

\begin{abstract}
The so-called pinball loss for estimating conditional quantiles
is a well-known tool in both statistics and machine learning.
So far, however, only little work has been done to quantify the efficiency of this
tool for nonparametric 
approaches.
We  fill this gap by establishing inequalities that describe
how close approximate pinball risk minimizers are to the corresponding conditional quantile.
These inequalities, which hold under mild assumptions on the data-generating distribution, are
then used to establish so-called variance bounds, which recently turned out to play an
important role in the statistical analysis of  (regularized) empirical risk minimization approaches.
Finally, we use both types of  inequalities
to establish an
oracle inequality for support vector machines that use the pinball loss.
The resulting learning rates are min--max  optimal
under some standard regularity assumptions on the   conditional quantile.
\end{abstract}


\end{frontmatter}

\section{Introduction}

Let $\mathrm{P}$ be a distribution on $X\times \mathbb{R}$, where $X$ is an arbitrary set equipped with a $\sigma$-algebra.
The goal of quantile regression is  to estimate the conditional quantile, that is, the set-valued function
\[
  F_{\tau,\mathrm{P}}^*(x) := \{ t\in \mathbb{R}\dvtx  \mathrm{P}((-\infty,t] | x) \geq \tau
   \mbox{ and } \mathrm{P}([t,\infty) | x)\geq
   1-\tau\}  , \qquad x\in  X ,
\]
where $\tau\in(0,1)$ is a fixed constant
specifying  the desired quantile level
and $\mathrm{P}( \cdot | x)$, $x\in X$,
is the regular conditional probability of $\mathrm{P}$.
 Throughout this paper, we assume
that $\mathrm{P}( \cdot | x)$ has its support in $[-1,1]$ for $\mathrm{P}_X$-almost all $x\in X$, where
$\mathrm{P}_X$ denotes the marginal distribution of $\mathrm{P}$ on $X$.
(By a simple scaling argument, all our results can be generalized to distributions living on
$X\times [-M,M]$ for some $M>0$. The uniform boundedness of the conditionals $\mathrm{P}( \cdot | x)$ is, however, crucial.)
Let us additionally assume for a moment
that
$F_{\tau,\mathrm{P}}^*(x)$ consists of singletons, that is, there exists an $f_{\tau,\mathrm{P}}^*\dvtx X\to \mathbb{R}$,
called the conditional $\tau$-quantile function,
such that $F_{\tau,\mathrm{P}}^*(x)= \{f_{\tau,\mathrm{P}}^*(x)\}$ for $\mathrm{P}_X$-almost all $x\in X$.
(Most of our main results do not require this assumption, but here, in the introduction, it makes the exposition more transparent.)
Then one approach to estimate the conditional $\tau$-quantile function is based
on the so-called \emph{$\tau$-pinball loss} $L\dvtx  Y   \times \mathbb{R}\to[0,\infty)$, which is defined by
\[
L(y,t):=
\cases{
(1-\tau) (t-y), &\quad \mbox{if }$ y < t$,\cr
\tau (y-t), &\quad \mbox{if }$ y \geq t$.
}
\]
With the help of this loss function we define
 the $L$-risk of a  function $f:X\to \mathbb{R}$
by
\[
\RP L f
:=
\mathbb{E}_{(x,y)\sim \mathrm{P}}L(y,f(x))
=
\int_{X\times Y} L (y,f(x))  \, \mathrm{d}\mathrm{P}(x,y)  .
\]
Recall  that $f_{\tau,\mathrm{P}}^*$ is up to  $\mathrm{P}_X$-zero   sets the \textit{only} function satisfying 
$\RP L {f_{\tau,\mathrm{P}}^*} = \inf  \RP L f =: \RPB L $, where the infimum is taken over \textit{all}
measurable functions   $f\dvtx X\to \mathbb{R}$.
Based on this observation, several estimators minimizing a (modified)  empirical $L$-risk were proposed
(see \cite{Koenker05} for a survey on both parametric and nonparametric methods) for situations where
$\mathrm{P}$ is unknown, but i.i.d.~samples
$D:= ((x_1,y_1),\dots,(x_n,y_n))  \in (X \times \mathbb{R})^n$ drawn from $\mathrm{P}$ are given.

Empirical methods estimating   quantile functions  with the help of the pinball loss
typically obtain  functions $f_D$ for which $\RP L {f_D}$ is close to $\RPB L$ with high probability.
In general, however, this only implies that $f_D$ is close to $f_{\tau,\mathrm{P}}^*$ in a very weak sense
(see \cite{Ste05b}, Remark 3.18) but recently, \cite {SteinwartChristmann2007b}, Theorem 2.5, established \textit{self-calibration
inequalities} of  the form
\begin{equation}\label{m1-ineq}
\snorm{f-f_{\tau,\mathrm{P}}^*}_{\Lx r{\mathrm{P}_X}} \leq
c_\mathrm{P}  \sqrt{ \RP {L} f - \RPB {L}  }  ,
\end{equation}
which hold under mild  assumptions on $\mathrm{P}$  described by the parameter $r\in (0,1]$.
The   first   goal of this paper is to generalize and to improve these inequalities.
Moreover, we will use these new self-calibration inequalities to establish
\textit{variance bounds} for the pinball risk, which in turn are known to improve the statistical
analysis of empirical risk minimization (ERM) approaches.

The second goal of this paper is to apply the self-calibration inequalities and the variance
bounds to support vector machines (SVMs) for quantile regression.
Recall that \cite{ScSmWiBa00aa,HwSh05a,TaLeSeSm06a} proposed an SVM  that finds a solution $f_{D,\lambda}\in H$ of
\begin{equation}\label{svm1}
  \arg\min_{f \in H}   \lambda \|f\|_H^2 +  \RT L f  ,
\end{equation}
where $\lambda>0$ is a regularization parameter, $H$ is a reproducing kernel Hilbert space (RKHS) over $X$, and
$\RT L f$ denotes the empirical risk of $f$, that is, $\RT L f:= \frac 1 n \sum_{i=1}^nL(y_i,f(x_i))$.
In \cite{ChMeSt09a} robustness properties and consistency for \textit{all} distributions $\mathrm{P}$ on $X\times \mathbb{R}$ were established
for this SVM, while
\cite{HwSh05a,TaLeSeSm06a} worked out how to  solve  this optimization problem with standard  techniques from  machine learning. Moreover,
\cite{TaLeSeSm06a} also provided an exhaustive empirical study, which shows
the excellent performance of this SVM.
We have recently established an oracle inequality for these SVMs in \cite{SteinwartChristmann2007b},
which was based on (\ref{m1-ineq}) and the resulting variance bounds.
In this paper, we improve this oracle inequality with the help of the new self-calibration
inequalities and variance bounds. It turns out that the resulting learning rates
are substantially faster than those of \cite{SteinwartChristmann2007b}.
Finally, we briefly discuss an adaptive parameter selection strategy.

The rest of this paper is organized as follows. In Section \ref{results}, we present both our new self-calibration inequality and
the new variance bound. We also introduce the assumptions on $\mathrm{P}$ that lead to these inequalities and discuss how these
inequalities improve our former results in \cite{SteinwartChristmann2007b}. In Section \ref{rates}, we use these new inequalities
to establish an oracle inequality for the SVM approach above. In addition,  we discuss the resulting learning rates and how these can be
achieved in an adaptive way.
Finally, all proofs
are contained in Section \ref{proofs}.

\section{Main results}\label{results}

In order to formulate the main results of this section, we need to
introduce some assumptions on the data-generating distribution $\mathrm{P}$.
To this end, let
 $\mathrm{Q}$ be a distribution on $\mathbb{R}$ and $\operatorname{supp} \mathrm{Q}$ be its support.
For $\tau\in (0,1)$, the \textit{$\tau$-quantile} of $\mathrm{Q}$ is the set
\[
  F_{\tau}^*(\mathrm{Q}) := \{ t\in \mathbb{R}:  \mathrm{Q}((-\infty,t]) \geq \tau  \mbox{ and }  \mathrm{Q}([t,\infty))\geq
   1-\tau\}  .
\]
It is well known that $F_{\tau}^*(\mathrm{Q})$ is a bounded  and closed interval.
We write
\[
t^*_{\mathrm{min}}(\mathrm{Q})  :=  \min F_{\tau}^*(\mathrm{Q}) \quad \mbox{and} \quad
t^*_{\mathrm{max}}(\mathrm{Q})  := \max F_{\tau}^*(\mathrm{Q})  ,
\]
which implies $F_{\tau}^*(\mathrm{Q}) = [t^*_{\mathrm{min}}(\mathrm{Q}), t^*_{\mathrm{max}}(\mathrm{Q})]$. Moreover,
it is easy to check that the interior of $F_{\tau}^*(\mathrm{Q})$ is a $\mathrm{Q}$-zero set,
that is, $\mathrm{Q}((t^*_{\mathrm{min}}(\mathrm{Q}), t^*_{\mathrm{max}}(\mathrm{Q}))) = 0$.
To avoid notational overload, we usually omit the
argument $\mathrm{Q}$ if the considered distribution is clearly
determined from the context.

\begin{definition}[(Quantiles of type $\bolds{q}$)]\label{distribut-type-Q}
A distribution $\mathrm{Q}$ with $\operatorname{supp} \mathrm{Q}\subset [-1,1]$
is said to have a $\tau$-quantile  of type $q\in (1,\infty)$ if
there exist
constants
$\alpha_\mathrm{Q}\in (0,2]$ and  $b_\mathrm{Q}>0$
such that
\begin{eqnarray}\label{q-type-1}
\mathrm{Q}\bigl( (t^*_{\mathrm{min}}-s,t^*_{\mathrm{min}})\bigr) &\geq & b_\mathrm{Q}
s^{q-1},
\\ \label{q-type-2}
\mathrm{Q}\bigl( (t^*_{\mathrm{max}},t^*_{\mathrm{max}}+s)\bigr) & \geq & b_\mathrm{Q}  s^{q-1}
\end{eqnarray}
for all $s\in [0,\alpha_\mathrm{Q}]$.
Moreover, $\mathrm{Q}$ has a $\tau$-quantile of type $q=1$,  if
$\mathrm{Q}(\{t^*_{\mathrm{min}}\})>0$ and \mbox{$\mathrm{Q}(\{t^*_{\mathrm{max}}\})>0$}. In this case, we define $\alpha_\mathrm{Q}:= 2$ and
\[
b_\mathrm{Q} :=
\cases{
         \min\{\mathrm{Q}(\{t^*_{\mathrm{min}}\}),\mathrm{Q}(\{t^*_{\mathrm{max}}\})\}, &
         \quad\mbox{if }$   t^*_{\mathrm{min}} \neq t^*_{\mathrm{max}}$,
         \cr
        \min\{\tau-\mathrm{Q}((-\infty,t^*_{\mathrm{min}})),\mathrm{Q}((-\infty,t^*_{\mathrm{max}}]) - \tau\},
         & \quad\mbox{if }$    t^*_{\mathrm{min}} = t^*_{\mathrm{max}}$,
}
\]
%
%
where we note that $b_\mathrm{Q}>0$ in both cases.
For all $q\geq 1$, we finally write $\gamma_\mathrm{Q} := b_\mathrm{Q} \alpha_\mathrm{Q}^{q-1}$.
\end{definition}

Since  $\tau$-quantiles of type $q$ are the central concept of this work, let us illustrate
this notion  by a few examples.
We begin with an example for which all quantiles are of type $2$.

\begin{example}
Let $\nu$ be a  distribution with $\operatorname{supp} \nu \subset [-1,1]$, $\mu$ be a distribution with
$\operatorname{supp} \mu \subset [-1,1]$ that has a density $h$
with respect to  the Lebesgue measure and $\mathrm{Q} := \alpha \nu + (1-\alpha)\mu$ for some  $\alpha\in [0,1)$.
If  $h$  is bounded away from $0$, that is, $h(y) \geq b$ for some $b>0$ and Lebesgue-almost all $y\in [-1,1]$, then
$\mathrm{Q}$ has a $\tau$-quantile of type $q=2$ for all $\tau\in (0,1)$ as simple integration shows. In this case, we set $b_\mathrm{Q}:= (1-\alpha)b$ and
$\alpha_\mathrm{Q}:=
\min\{1+t^*_{\mathrm{min}},1-t^*_{\mathrm{max}}\}$.
\end{example}


\begin{example}
Again, let $\nu$ be a distribution with $\operatorname{supp} \nu \subset [-1,1]$, $\mu$ be a distribution with
$\operatorname{supp} \mu \subset [-1,1]$ that has a
Lebesgue density $h$,
 and $\mathrm{Q} := \alpha \nu + (1-\alpha)\mu$ for some  $\alpha\in [0,1)$.
If, for a fixed $\tau\in (0,1)$, there exist constants $b>0$ and $p>-1$ such that
\begin{eqnarray*}
 h(y) & \geq & b  \bigl(t^*_{\mathrm{min}}(\mathrm{Q}) - y \bigr)^p   ,\qquad  y\in [-1,t^*_{\mathrm{min}}(\mathrm{Q})],
 \\
 h(y) & \geq & b  \bigl(y-t^*_{\mathrm{max}}(\mathrm{Q})  \bigr)^p  , \qquad  y\in
 [t^*_{\mathrm{max}}(\mathrm{Q}),1].
\end{eqnarray*}
Lebesgue-almost surely, then simple integration shows that $\mathrm{Q}$ has a $\tau$-quantile of type $q=2+p$ and
we may set $b_\mathrm{Q}:=(1-\alpha) b/(1+p)$ and $\alpha_\mathrm{Q}:= \min\{1+t^*_{\mathrm{min}}(\mathrm{Q}),1-t^*_{\mathrm{max}}(\mathrm{Q})\}$.
\end{example}

\begin{example}
Let $\nu$ be a distribution with $\operatorname{supp} \nu \subset [-1,1]$ and
$\mathrm{Q} := \alpha \nu + (1-\alpha) \delta_{t^*}$ for some  $\alpha\in [0,1)$, where $\delta_{t^*}$ denotes the Dirac measure at $t^*\in (0,1)$.
If $\nu(\{t^*\}) = 0$, we then have $\mathrm{Q}((-\infty,t^*)) = \alpha \nu  ((-\infty,t^*))$ and
$\mathrm{Q}((-\infty,t^*]) = \alpha \nu  ((-\infty,t^*)) + 1-\alpha$, and hence $\{t^*\}$ is a $\tau$-quantile of type $q=1$
for all $\tau$ satisfying $\alpha \nu  ((-\infty,t^*)) < \tau < \alpha \nu  ((-\infty,t^*)) + 1-\alpha$.
\end{example}

\begin{example}
Let $\nu$ be a distribution with $\operatorname{supp} \nu \subset [-1,1]$ and
$\mathrm{Q} := (1-\alpha-\beta) \nu + \alpha \delta_{t_{\mathrm{min}}} + \beta \delta_{t_{\mathrm{max}}}$ for some
$\alpha,\beta\in (0,1]$ with $\alpha+\beta\leq 1$.
If $\nu([t_{\mathrm{min}},t_{\mathrm{max}}]) = 0$, we  have
$\mathrm{Q}((-\infty,t_{\mathrm{min}}]) = (1-\alpha-\beta) \nu  ((-\infty,t_{\mathrm{min}}]) + \alpha$ and
$\mathrm{Q}([t_{\mathrm{max}},\infty)) = (1-\alpha-\beta) (1-\nu  ((-\infty,t_{\mathrm{min}}]))  + \beta$. Consequently,
 $[t_{\mathrm{min}},t_{\mathrm{max}}]$
is the $\tau:= (1-\alpha-\beta) \nu  ((-\infty,t^*])  +\alpha$
quantile of $\mathrm{Q}$ and this quantile is of type $q=1$.
\end{example}

As outlined in the introduction,
we are not interested in a single distribution $\mathrm{Q}$ on $\mathbb{R}$
but in distributions $\mathrm{P}$ on $X\times \mathbb{R}$. The following definition
extends the previous definition to such $\mathrm{P}$.

\begin{definition}[(Quantiles of  $\bolds{p}$-average type $\bolds{q}$)]
Let $p\in (0,\infty]$,  $q\in [1,\infty)$, and
$\mathrm{P}$ be a distribution on $X\times \mathbb{R}$ with
 $\operatorname{supp} \mathrm{P}(\cdot |x)\subset [-1,1]$ for $\mathrm{P}_X$-almost all $x\in X$.
Then $\mathrm{P}$ is
%
said to have
a $\tau$-quantile of
$p$-average type $q$, if
$\mathrm{P}(\cdot |x)$ 
has a $\tau$-quantile of type $q$ for $\mathrm{P}_X$-almost all $x\in X$, and
the function $\gamma\dvtx X\to [0,\infty]$ defined, for $\mathrm{P}_X$-almost all $x\in X$, by
\[
\gamma(x):= \gamma_{\mathrm{P}(\cdot|x)}  ,
\]
where $\gamma_{\mathrm{P}(\cdot|x)} = b_{\mathrm{P}(\cdot|x)} \alpha_{\mathrm{P}(\cdot|x)}^{q-1}$
is defined in Definition \ref{distribut-type-Q}, satisfies
$\gamma^{-1} \in \Lx p {\mathrm{P}_X}$.
\end{definition}

To  establish the announced self-calibration inequality,
we finally need the distance
\[
\operatorname{dist}(t,A) := \inf_{s\in A} |t-s|
\]
between an element  $t\in \mathbb{R}$  and an $A\subset \mathbb{R}$.
Moreover, $\operatorname{dist}(f,F^*_{\tau,\mathrm{P}})$ denotes the function $x\mapsto \operatorname{dist}(f(x),F^*_{\tau,\mathrm{P}}(x))$.
With these preparations
the self-calibration inequality reads as follows.

\begin{theorem}\label{main1}
Let $L$ be the $\tau$-pinball loss,
 $p\in (0,\infty]$ and $q\in [1,\infty)$ be real numbers, and $r:= \frac {pq} {p+1}$. Moreover, let $\mathrm{P}$
be a distribution that has a $\tau$-quantile of   $p$-average type $q\in [1,\infty)$.
Then, for all
$f\dvtx X\to [-1,1]$,
we have
\[
\snorm{\operatorname{dist}(f,F^*_{\tau,\mathrm{P}})}_{\Lx r{\mathrm{P}_X}}
\leq
2^{1-1/q}q^{1/q} \snorm{\gamma^{-1}}_{L_p(\mathrm{P}_X)}^{1/q}  \bigl(\RP {L} f - \RPB {L}  \bigr)^{1/q}  .
\]
\end{theorem}

Let us briefly compare the self-calibration inequality above with the one
established in  \cite{SteinwartChristmann2007b}. To this end, we can solely
focus on the case $q=2$, since this was the only case considered in
\cite{SteinwartChristmann2007b}. For the same reason, we can restrict our considerations
to distributions $\mathrm{P}$ that have a unique conditional $\tau$-quantile $f_{\tau,\mathrm{P}}^*(x)$
for $\mathrm{P}_X$-almost all $x\in X$.
Then Theorem \ref{main1} yields
\[
\snorm{f-f_{\tau,\mathrm{P}}^*}_{\Lx r{\mathrm{P}_X}} \leq 2
\snorm{\gamma^{-1}}_{L_p(\mathrm{P}_X)}^{1/2}  \bigl(\RP {L} f - \RPB {L}  \bigr)^{1/2}
\]
for $r:= \frac {2p} {p+1}$. On the other hand,   it was shown in
\cite{SteinwartChristmann2007b}, Theorem 2.5, that
\[
\snorm{f-f_{\tau,\mathrm{P}}^*}_{\Lx {r/2}{\mathrm{P}_X}} \leq \sqrt 2
\snorm{\gamma^{-1}}_{L_p(\mathrm{P}_X)}^{1/2}  \bigl(\RP {L} f - \RPB {L}  \bigr)^{1/2}
\]
under the \textit{additional} assumption that the conditional widths $\alpha_{\mathrm{P}(\cdot|x)}$
considered in Definition \ref{distribut-type-Q} are \textit{independent} of $x$.
Consequently,  our new self-calibration inequality is
more general and, modulo the constant $\sqrt 2$,  also sharper.

It is well known that self-calibration inequalities for Lipschitz continuous losses
lead to variance bounds, which in turn are important for the
statistical analysis of ERM approaches; see
\cite{MaTs99a,Massart00b,Mendelson01a,Mendelson01b,Tsybakov04a,BaBoMe05a,BaJoMc03a}.
For the pinball loss, we obtain the following variance bound.


\begin{theorem}\label{main2}
Let $L$ be the $\tau$-pinball loss,
 $p\in (0,\infty]$ and $q\in [1,\infty)$ be real numbers, and
\[
 \vartheta:=\min\biggl\{\frac 2 q ,   \frac {p} {p+1}\biggr\}  .
\]
Let $\mathrm{P}$
be a distribution that has a $\tau$-quantile of   $p$-average type $q$.
Then, for all
$f\dvtx X\to [-1,1]$, there exists an $f^*_{\tau,\mathrm{P}}\dvtx X\to [-1,1]$
with $f^*_{\tau,\mathrm{P}}(x) \in  F^*_{\tau,\mathrm{P}}(x)$ for $\mathrm{P}_X$-almost all $x\in X$ such that
\[
\mathbb{E}_\mathrm{P} (  L\circ f - L\circ \fpb     )^2
\leq
2^{2-\vartheta} q^{\vartheta}
\snorm{\gamma^{-1}}_{L_p(\mathrm{P}_X)}^{\vartheta}  \bigl(\RP L f - \RPB L \bigr)^{\vartheta}  ,
\]
where we used the shorthand $L\circ f$ for the function $(x,y)\mapsto L(y,f(x))$.
\end{theorem}

Again, it is straightforward to show that the variance bound above is both more general and stronger than the variance
bound established in \cite {SteinwartChristmann2007b}, Theorem 2.6.

\section{An application to support vector machines}\label{rates}

The goal of this section is to establish an oracle inequality for
the SVM defined in (\ref{svm1}).
The use of this oracle inequality is then illustrated by some learning rates
we derive from it.

Let us begin by recalling some RKHS theory (see, e.g., \cite {SteinwartChristmann2008a}, Chapter 4, for a more detailed account).
 To this end, let $k\dvtx X\times X\to \mathbb{R}$
be a   measurable kernel, that is, a measurable function that is  symmetric and positive definite.
Then the associated RKHS $H$ consists of measurable functions. Let us additionally assume that $k$ is bounded with
$\inorm{k}:=\sup_{x\in X}\sqrt{k(x,x)}\leq 1$,  which in turn implies that $H$ consists of bounded functions and
$\inorm{f} \leq  \hnorm{f}$ for all $f\in H$.

Suppose now that we have a distribution $\mathrm{P}$ on $X\times Y$.
To describe the approximation error of SVMs we
use the \emph{approximation  error function}
\[
A(\lambda) := \inf_{f\in H} \lambda \snorm f_H^2 + \RP L f - \RPB L  , \qquad \lambda>0,
\]
where $L$ is the $\tau$-pinball loss.
Recall that \cite {SteinwartChristmann2008a}, Lemma 5.15 and Theorem 5.31,  showed that
$\lim_{\lambda\to 0} A(\lambda) = 0$, if the RKHS $H$ is dense in $\Lx 1 {\mathrm{P}_X}$
and the speed of this convergence describes how well $H$ approximates the Bayes $L$-risk $\RPB L$.
In particular, \cite {SteinwartChristmann2008a}, Corollary 5.18, shows that
$A(\lambda) \leq c \lambda$ for some constant $c>0$ and all $\lambda>0$ if and only if there exists an
$f\in H$ such that $f(x) \in  F^*_{\tau,\mathrm{P}}(x)$ for $\mathrm{P}_X$-almost all $x\in X$.

We further  need the
integral operator $T_k\dvtx\Lx 2 {\mathrm{P}_X}\to \Lx 2 {\mathrm{P}_X}$   defined by
\[
T_k f (\cdot):= \int_X k(x,\cdot) f(x)  \, \mathrm{d}\mathrm{P}_X(x)   ,  \qquad f\in \Lx 2 {\mathrm{P}_X}.
\]
It is well known that $T_k$ is self-adjoint and nuclear; see, for example, \cite {SteinwartChristmann2008a}, Theorem 4.27.
Consequently, it has at most countably many eigenvalues (including geometric multiplicities),
which are all non-negative and summable.
Let us order these eigenvalues $\lambda_i(T_k)$.
Moreover, if we only have finitely many eigenvalues, we extend this finite sequence
by zeros. As a result, we  can always deal with a decreasing, non-negative sequence
$\lambda_1(T_k)\geq  \lambda_2(T_k)\geq \cdots$, which satisfies $\sum_{i=1}^\infty \lambda_i(T_k) < \infty$.
The finiteness of this sum can already be used to establish oracle inequalities; see~\cite {SteinwartChristmann2008a}, Theorem 7.22.
But in the following we assume that the eigenvalues converge even faster to zero, since
{(a)} this case  is satisfied for many RKHSs and (b) it leads to better
oracle inequalities.
To be more precise, we assume that there exist constants $a\geq 1$ and $\varrho\in (0,1)$ such that
\begin{equation}\label{eigenvalues-assump}
\lambda_i (T_k) \leq a  i^{-1/\varrho}  , \qquad  i\geq 1.
\end{equation}
Recall that  (\ref{eigenvalues-assump}) was first used in \cite{BlBoMa04a}
to establish an oracle inequality for SVMs using the hinge loss, while \cite{CaDe07a,MeNe08a,StHuSc09b}
consider (\ref{eigenvalues-assump}) for SVMs using the least-squares loss.
Furthermore, one can show (see \cite{Ste08a}) that (\ref{eigenvalues-assump}) is equivalent (modulo a  constant only depending on $\varrho$) to
\begin{equation}\label{EN}
e_i\bigl(\operatorname{id}:H\to \Lx 2 {\mathrm{P}_X}\bigr) \leq \sqrt a i^{-1/(2\varrho)}  , \qquad  i\geq 1  ,
\end{equation}
where $e_i(\operatorname{id}\dvtx H\to \Lx 2 {\mathrm{P}_X})$ denotes the $i$th (dyadic) entropy number \cite{CaSt90} of the inclusion map from $H$ into $\Lx 2 {\mathrm{P}_X}$.
In addition, \cite{Ste08a} shows that (\ref{EN}) implies a bound on expectations of
random entropy numbers, which in turn are used in \cite {SteinwartChristmann2008a}, Chapter 7.4, to establish general oracle inequalities for SVMs.
On the other hand, (\ref{EN}) has been extensively studied in the literature. For example,
 for $m$-times differentiable kernels on Euclidean balls $X$ of $\mathbb{R}^d$,
it is known that (\ref{EN}) holds for $\varrho:= \frac d{2m}$.  We refer to~\cite {CuZh07}, Chapter 5, and
 \cite {SteinwartChristmann2008a}, Theorem 6.26, for a precise statement.
 Analogously, if  $m>d/2$ is some integer, then
 the Sobolev space $H:=W^m(X)$ is an RKHS that satisfies (\ref{EN}) for $\varrho:= \frac {d}{2m}$, and
this estimate is also asymptotically sharp; see  \cite{BiSo67a,EdTr96}.

We finally
need the clipping operation 
defined by
\[
\clippt := \max\{ -1, \min\{ 1, t\}\}
\]
%
%
for all $t\in \mathbb{R}$.
We can now state
the following oracle inequality for SVMs using the pinball loss.


\begin{theorem}\label{oracle-general}
Let $L$ be the $\tau$-pinball loss and
 $\mathrm{P}$ be a  distribution on $X\times \mathbb{R}$ with $\operatorname{supp} \mathrm{P}(\cdot |x)\subset [-1,1]$ for
$\mathrm{P}_X$-almost all $x\in X$.
Assume that there exists a function
$f_{\tau,\mathrm{P}}^*\dvtx X\to \mathbb{R}$ with
$f_{\tau,\mathrm{P}}^*(x) \in  F^*_{\tau,\mathrm{P}}(x)$ for $\mathrm{P}_X$-almost all $x\in X$ and
 constants $V \geq 2^{2-\vartheta}$ and  $\vartheta\in[0,1]$
such that
\begin{equation}\label{var-bound}
\mathbb{E}_\mathrm{P} (  L\circ f - L\circ \fpb     )^2
\leq
V   \bigl(\RP L f - \RPB L  \bigr)^\vartheta
\end{equation}
for all  $f\dvtx X\to [-1,1]$.
Moreover, let $H$ be a separable RKHS over $X$ with a bounded measurable kernel satisfying $\inorm k\leq 1$.
In addition, assume that (\ref{eigenvalues-assump}) is satisfied for some $a\geq 1$ and $\varrho\in (0,1)$.
Then there exists a constant $K$ depending only on $\varrho$, $V$, and $\vartheta$
such that, for all $\varsigma \geq 1$, $n
\geq 1$ and $\lambda >0 $,
we have
with probability $\mathrm{P}^n$ not less than $1-3e^{-\varsigma}$ that
\[
 \RP L \fTc - \RPB L
 \leq
9 A(\lambda) + 30 \sqrt{\frac {A(\lambda)} \lambda}\frac {\varsigma}{n}
+ K \biggl(\frac {a^{\varrho}}{\lambda ^\varrho n}  \biggr)^{ 1/( {2-\varrho-\vartheta+\vartheta \varrho})}
+ 3 \biggl(\frac{72 V\varsigma}n\biggr)^{ 1/({2-\vartheta})}  .
\]
\end{theorem}

Let us now discuss the learning rates obtained from this oracle inequality.
%
To this end, we  assume in the following  that there exist constants $c>0$ and $\beta\in (0,1]$ such that
\begin{equation}\label{A2}
A(\lambda)\leq c\lambda^\beta  ,  \qquad \lambda>0.
\end{equation}
Recall from \cite {SteinwartChristmann2008a}, Corollary 5.18, that, for $\beta=1$,
this assumption holds if and only if there exists a $\tau$-quantile function
$f_{\tau,\mathrm{P}}^*$ with $f_{\tau,\mathrm{P}}^* \in H$. Moreover, for $\beta<1$, there is a tight relationship between (\ref{A2})
and the behavior of the  approximation error of the balls $\lambda^{-1}B_H$; see \cite{SteinwartChristmann2008a}, Theorem~5.25.
In addition, one can show (see \cite{SteinwartChristmann2008a}, Chapter 5.6) that if $f_{\tau,\mathrm{P}}^*$ is contained in
the real interpolation space $(\Lx 1 {\mathrm{P}_X},H)_{\vartheta,\infty}$, see  \cite{BeSh88}, then (\ref{A2}) is satisfied for
$\beta:=\vartheta /(2-\vartheta)$. For example, if $H:=W^m(X)$ is a Sobolev space over a Euclidean ball $X\subset \mathbb{R}^d$ of
order $m>d/2$ and $\mathrm{P}_X$
has a Lebesgue density that is bounded away from $0$ and $\infty$, then $f_{\tau,\mathrm{P}}^* \in W^s(X)$ for some $s\in (d/2,m]$ implies
(\ref{A2}) for $\beta:=   s /({2m-s})$.

Now assume that (\ref{A2}) holds. We further assume that $\lambda$ is determined by $\lambda_n = n^{-\gamma/\beta}$, where
\begin{equation}\label{def-gamma}
\gamma:= \min\biggl\{\frac{\beta}{\beta(2-\vartheta+\varrho\vartheta-\varrho)+\varrho},  \frac{2\beta}{\beta+1}\biggr\}  .
\end{equation}
Then Theorem \ref{oracle-general} shows that
$\RP L {\fTcn}$ converges to  $\RPB L$ with rate $n^{-\gamma}$;
see  \cite{SteinwartChristmann2008a}, Lemma~A.1.7, for
calculating the value of $\gamma$.
%
Note that this choice of $\lambda$ yields the best learning rates from Theorem \ref{oracle-general}. Unfortunately,
however, this choice  requires knowledge of the usually unknown parameters $\beta$, $\vartheta$ and $\varrho$.
To address this issue, let us consider the following scheme that is close to approaches taken in practice (see
\cite{Rosset07a} for a similar technique that has a fast implementation based on regularization paths).

\begin{definition}\label{conc-bas:tv-svm}
Let $H$ be an RKHS over $X$
and  $\Lambda :=(\Lambda_n)$ be a sequence of finite subsets $\Lambda_n \subset (0,1]$.
Given a data set
 $D:=((x_{1},y_{1}),\dots,(x_n,y_n))\in (X\times \mathbb{R})^n$, we define
 \begin{eqnarray*}
 D_1 &:=& ((x_{1},y_{1}),\dots,(x_{m},y_{m}))  ,
  \\
D_2  &:=& ((x_{m+1},y_{m+1}),\dots,(x_{n},y_{n}))  ,
\end{eqnarray*}
where $m:= \lfloor n/2\rfloor + 1$ and $n\geq 3$. Then we use
 $D_1$ to compute the SVM decision functions
\[
f_{D_1,\lambda} := \arg\min_{f\in H} \lambda \snorm f_H^2 + \Rx L{\mathrm{D}_1} f  , \qquad  \lambda \in \Lambda_n,
\]
and  $D_2$ to determine $\lambda$ by  choosing a $\lambda_{D_2}\in \Lambda_n$ such  that\vspace*{-1ex}
\[
\Rx L {\mathrm{D}_2} {\clippfo_{D_1,\lambda_{D_2}}} = \min_{\lambda\in \Lambda_n} \Rx L {\mathrm{D}_2} {\clippfo_{D_1,\lambda}}   .\vspace*{-1ex}
\]
In the following, we call this learning method, which produces
 the  decision functions  $\clippfo_{D_1,\lambda_{D_2}}$,   a training validation SVM  with respect to $\Lambda$.
\end{definition}

Training validation SVMs have been extensively studied in \cite{SteinwartChristmann2008a}, Chapter 7.4.
In particular, \cite{SteinwartChristmann2008a}, Theorem 7.24, gives the following result
that shows that  the
learning rate $n^{-\gamma}$ can be achieved without knowing of the existence of the
parameters $\beta$, $\vartheta$ and $\varrho$ or their particular values.

\begin{theorem}\label{tv-svm-rate}
Let $(\Lambda_n)$ be a sequence of
$n^{-2}$-nets $\Lambda_n$ of $(0,1]$ such that the cardinality $|\Lambda_n|$ of $\Lambda_n$ grows polynomially in $n$.
Furthermore, consider the situation of Theorem \ref{oracle-general} and assume that (\ref{A2}) is satisfied for some $\beta\in (0,1]$.
Then the   training validation SVM  with respect to $\Lambda :=(\Lambda_n)$ learns with rate $ n^{- \gamma}$, where $\gamma$ is defined by
(\ref{def-gamma}).
\end{theorem}

Let us now consider how these learning rates in terms of risks translate into
rates for
\begin{equation}\label{rnorm}
\snorm{\fTcn - f_{\tau,\mathrm{P}}^*}_{\Lx r {\mathrm{P}_X}} \to 0  .
\end{equation}
To this end, we assume that $\mathrm{P}$ has a $\tau$-quantile of $p$-average type $q$,
where we additionally assume for the sake of simplicity that
 $r:=\frac{pq}{p+1} \le 2$.
Note that the latter is satisfied for all $p$ if $q\leq 2$, that is, if all conditional distributions are concentrated
around the quantile at least as much as the uniform distribution; see the discussion following Definition  \ref{distribut-type-Q}.
We  further  assume that the conditional quantiles $F^*_{\tau,\mathrm{P}}(x)$ are singletons for
$\mathrm{P}_X$-almost all $x\in X$. Then Theorem \ref{main2} provides a variance bound of the form (\ref{var-bound})
for $\vartheta := p/(p+1)$, and hence $\gamma$ defined in (\ref{def-gamma}) becomes
\[
\gamma= \min\biggl\{\frac{\beta(p+1)}{\beta(2+p - \varrho)+\varrho (p+1) },  \frac{2\beta}{\beta+1}\biggr\}  .
\]
By Theorem \ref{main1} we consequently see that  (\ref{rnorm})
converges with rate $n^{-\gamma/q}$, where $r:=pq/\break(p+1)$.
To illustrate this learning rate, let us assume that we have picked an RKHS $H$ with
$f_{\tau,\mathrm{P}}^*\in H$. Then we have $\beta=1$, and hence
it is easy to check that the latter learning rate reduces~to
\[
 n^{-{(p+1)/(q(2+p+\varrho p))}}.
\]
For the sake of simplicity, let us further assume  that the conditional distributions
do not change too much in the sense that $p=\infty$.
Then we have $r=q$, and hence
\begin{equation}\label{excess}
\int_X |\fTcn - f_{\tau,\mathrm{P}}^* |^q\,   \mathrm{d}\mathrm{P}_X
\end{equation}
converges to zero with rate $n^{-1/(1+\varrho)}$.
The latter shows that the value of $q$ does not change the learning rate for (\ref{excess}), but only the exponent in (\ref{excess}).
Now note that by our assumption on $\mathrm{P}$ and the definition of the clipping operation we have\vspace*{-1ex}
\[
\inorm{\fTcn - f_{\tau,\mathrm{P}}^*}\leq 2  ,\vspace*{-1ex}
\]
and consequently small values of $q$ emphasize the discrepancy of $\fTcn$ to $f_{\tau,\mathrm{P}}^*$
more than large values of $q$ do.
In this sense, a stronger average concentration around the quantile  is helpful for the learning process.

Let us now have a closer look at   the special case $q=2$, which is
probably the most interesting case for applications.
Then  we have the learning rate
$n^{-1/(2(1+\varrho))}$ for\vspace*{-1ex}
\[
\snorm{\fTcn - f_{\tau,\mathrm{P}}^*}_{\Lx 2 {\mathrm{P}_X}}  .
\]
Now recall that  the conditional median equals the conditional mean
for \emph{symmetric} conditional distributions   $\mathrm{P}(\cdot|x)$.  Moreover, if $H$ is  a Sobolev space
$W^m(X)$, where $m>d/2$ denotes the smoothness index and $X$ is a Euclidean ball in $\mathbb{R}^d$, then $H$ consists of continuous functions,
and
\cite{EdTr96} shows that $H$ satisfies (\ref{eigenvalues-assump}) for $\varrho:= d/(2m)$. 
Consequently,
 we see that in this case the latter convergence rate is optimal in a min--max sense \cite{YaBa99a,Temlyakov04a} if $\mathrm{P}_X$ is the uniform distribution.
Finally, recall that in the case $\beta=1$, $q=2$ and $p=\infty$ discussed so far,   the results derived in \cite{SteinwartChristmann2007b}   only yield a learning rate of $n^{-1/(3(1+\varrho))}$ for \vspace*{-1ex}
\[
\snorm{\fTcn - f_{\tau,\mathrm{P}}^*}_{\Lx 1 {\mathrm{P}_X}}  .
\]
In other words,  the earlier rates from \cite{SteinwartChristmann2007b}  are
not only worse by a factor of $3/2$ in the exponent but also are stated in terms of the weaker $\Lx 1 {\mathrm{P}_X}$-norm.
In addition, \cite{SteinwartChristmann2007b}  only considered the case $q=2$, and
hence we see that our new results are also more general.

\section{Proofs}\label{proofs}

Since the proofs of Theorems \ref{main1} and \ref{main2} use some notation developed in
\cite{Ste05b}
and \cite{SteinwartChristmann2008a}, Chapter~3, let us begin
by recalling these. To this end, let $L$ be the $\tau$-pinball loss for some fixed $\tau\in (0,1)$ and $\mathrm{Q}$ be
a distribution on $\mathbb{R}$ with $\operatorname{supp} \mathrm{Q}\subset [-1,1]$. Then \cite{Ste05b,SteinwartChristmann2008a} defined the
 \textit{inner $L$-risks}
by
\[
\CQ{L}{t} := \int_{Y} L(y,t)  \,\mathrm{d}\mathrm{Q}(y)  ,  \qquad  t\in\mathbb{R},
\]
and the  \textit{minimal inner $L$-risk} was denoted by
$\CQB L  := \inf_{t\in \mathbb{R}}  \CQ L  t$.
Moreover, we write
$\FQ L {0^{+}} = \{t\in \mathbb{R}\dvtx \CQ L t = \CQB L\}$ for the set of exact minimizers.

Our first goal is
to compute the excess inner risks and the set of  exact minimizers
for the pinball loss.
To this end recall
that (see \cite{Bauer01}, Theorem 23.8), given a distribution $\mathrm{Q}$ on $\mathbb{R}$ and a
 measurable  function
$g\dvtx X\to [0,\infty),$ we have
\begin{equation}\label{bauer-h1}
  \int_\mathbb{R} g \,\mathrm{d}\mathrm{Q} = \int_0^\infty \mathrm{Q}(\{g\geq s\}) \,  \mathrm{d}s  .
\end{equation}
With   these preparations we can now show the following generalization of \cite{SteinwartChristmann2008a}, Proposition 3.9.

\begin{proposition}\label{loss:pin-ball-more}
Let $L$ be the $\tau$-pinball loss and
$\mathrm{Q}$ be a distribution on $\mathbb{R}$ with $\CQB L<\infty$. 
Then there exist $q_+, q_- \in[0,1]$ with $q_++q_-=\mathrm{Q}([t^*_{\mathrm{min}},t^*_{\mathrm{max}}])$, and, for all $t\geq 0$, we have
\begin{eqnarray}\label{loss:pin-ball-more-a1}
\CQ L {t^*_{\mathrm{max}}+t} - \CQB L & = &
 tq_++\int_0^t \mathrm{Q}\bigl( (t^*_{\mathrm{max}},t^*_{\mathrm{max}}+s)\bigr) \,  \mathrm{d}s   ,
\\ \label{loss:pin-ball-more-a2}
\CQ L {t^*_{\mathrm{min}}-t} -  \CQB L & = &
{t}   q_-+\int_0^{t}    \mathrm{Q}\bigl( (t^*_{\mathrm{min}}-s,t^*_{\mathrm{min}})\bigr)\,  \mathrm{ d}s   .
\end{eqnarray}
Moreover, if $t^*_{\mathrm{min}} \neq t^*_{\mathrm{max}}$, then we have  $q_- = \mathrm{Q}(\{t^*_{\mathrm{min}}\})$ and $q_+ = \mathrm{Q}(\{t^*_{\mathrm{max}}\})$.
Finally, $\FQ L {0^+}$ equals the $\tau$-quantile, that is,  $\FQ L {0^+} = [t^*_{\mathrm{min}},t^*_{\mathrm{max}}]$.
\end{proposition}

\begin{pf}
Obviously, we have
$\mathrm{Q}  ( (-\infty,t^*_{\mathrm{max}}]) + \mathrm{Q}  ( [t^*_{\mathrm{max}},\infty)) = 1 + \mathrm{Q}(\{t^*_{\mathrm{max}}\})$, and hence we obtain
$\tau\leq \mathrm{Q} ( (-\infty,t^*_{\mathrm{max}}] ) \leq \tau + \mathrm{Q}(\{t^*_{\mathrm{max}}\})$. In other words, there
exists a $q_+\in[0,1]$   satisfying $0\leq q_+\leq \mathrm{Q}(\{t^*_{\mathrm{max}}\})$ and
\begin{equation}\label{loss:pin-ball-more-h1}
\mathrm{Q} ( (-\infty,t^*_{\mathrm{max}}] )  = \tau + q_+  .
\end{equation}
Let us consider the distribution $\tilde \mathrm{Q}$ defined by $\tilde \mathrm{Q} (A) := \mathrm{Q}(t^*_{\mathrm{max}}+A)$
for all measurable \mbox{$A\subset \mathbb{R}$}.
Then it is not hard to see that $t^*_{\mathrm{max}}(\tilde \mathrm{Q}) = 0$. Moreover, we obviously have
$\CQ L {t^*_{\mathrm{max}}+t} = \Cxx L{\tilde \mathrm{Q}} t$ for all $t\in \mathbb{R}$.
Let us now compute the inner risks of $L$  with respect to $\tilde \mathrm{Q}$.
To this end,  we fix a $t\geq 0$.
Then we have
\[
\int_{y<t}(y-t)   \,\mathrm{d}\tilde \mathrm{Q}(y)  = \int_{y<0} y  \, \mathrm{d}\tilde \mathrm{Q}(y) - t\tilde \mathrm{Q}((-\infty,t))
+ \int_{0\leq y<t} y  \, \mathrm{d}\tilde \mathrm{Q}(y)
\]
and
\[
\int_{y\geq t}(y-t)  \, \mathrm{d}\tilde \mathrm{Q}(y) =  \int_{y\geq 0}  y \,  \mathrm{d}\tilde \mathrm{Q}(y)
- t\tilde \mathrm{Q}([t,\infty)) - \int_{0\leq y<t} y   \,\mathrm{d}\tilde \mathrm{Q}(y)
\]
and hence we obtain
\begin{eqnarray*}
 \Cxx L{\tilde \mathrm{Q}}  {t} & = &
  (\tau-1)\int_{y<t}(y-t)  \, \mathrm{d}\tilde \mathrm{Q}(y) + \tau \int_{y\geq t}(y-t)\,   \mathrm{d}\tilde \mathrm{Q}(y) \nonumber
  \\
  & = &
    \Cxx L{\tilde \mathrm{Q}}  0 - \tau t + t\tilde \mathrm{Q}((-\infty,0)) + t\tilde \mathrm{Q}([0,t))
    - \int_{0\leq y<t} y  \, \mathrm{d}\tilde \mathrm{Q}(y)  .
   \nonumber
\end{eqnarray*}
Moreover, using (\ref{bauer-h1}) we find
\[
t\tilde \mathrm{Q}([0,t)) - \int_{0\leq y<t} y  \, \mathrm{d}\tilde \mathrm{Q}(y)
 =
 \int_0^t\tilde \mathrm{Q}([0,t))\,\mathrm{d}s  - \int_0^t \tilde \mathrm{Q}([s,t)) \, \mathrm{d}s
 = t \tilde \mathrm{Q}(\{0\}) +  \int_0^t\tilde \mathrm{Q}((0,s))\,\mathrm{d}s  ,
\]
and since (\ref{loss:pin-ball-more-h1}) implies $\tilde \mathrm{Q}((-\infty,0)) + \tilde \mathrm{Q}(\{0\})
=  \tilde \mathrm{Q}((-\infty,0]) = \tau + q_+$, we thus obtain
\begin{equation}\label{loss:pin-ball-more-a1-h1}
\CQ L {t^*_{\mathrm{max}}+t}   =    \CQ L t^*_{\mathrm{max}}+
 tq_++\int_0^t \mathrm{Q}\bigl( (t^*_{\mathrm{max}},t^*_{\mathrm{max}}+s)\bigr)   \,\mathrm{d}s  .
\end{equation}
By considering the pinball loss with parameter $1-\tau$ and the distribution
$\bar \mathrm{Q}$ defined by
$\bar \mathrm{Q}(A) := \mathrm{Q}(-t^*_{\mathrm{min}}-A)$, $A\subset \mathbb{R}$ measurable, we further see that (\ref{loss:pin-ball-more-a1-h1}) implies
\begin{equation}\label{loss:pin-ball-more-a1-h2}
\CQ L {t^*_{\mathrm{min}}-t}   =    \CQ L t^*_{\mathrm{min}}+
 {t}   q_-+\int_0^{t}    \mathrm{Q}\bigl( (t^*_{\mathrm{min}}-s,t^*_{\mathrm{min}})\bigr)\,   \mathrm{d}s   ,  \qquad t\geq 0,
\end{equation}
where $q_-$ satisfies $0\leq q_-\leq \mathrm{Q}(\{t^*_{\mathrm{min}}\})$
and $\mathrm{Q}([t^*_{\mathrm{min}},\infty)) = 1-\tau+q_-$. By
(\ref{loss:pin-ball-more-h1}) we then find
$q_++q_- = \mathrm{Q}([t^*_{\mathrm{min}},t^*_{\mathrm{max}}])$.
Moreover, if $t^*_{\mathrm{min}} \neq t^*_{\mathrm{max}}$, the fact  $\mathrm{Q}((t^*_{\mathrm{min}}, t^*_{\mathrm{max}})) = 0$ yields
\[
q_++q_-=\mathrm{Q}([t^*_{\mathrm{min}}, t^*_{\mathrm{max}}]) = \mathrm{Q}(\{t^*_{\mathrm{min}}\}) + \mathrm{Q}(\{t^*_{\mathrm{max}}\})  .
\]
Using the earlier established
$q_+\leq \mathrm{Q}(\{t^*_{\mathrm{max}}\})$ and $q_-\leq \mathrm{Q}(\{t^*_{\mathrm{min}}\})$, we then find both
$q_- = \mathrm{Q}(\{t^*_{\mathrm{min}}\})$ and $q_+ = \mathrm{Q}(\{t^*_{\mathrm{max}}\})$.

To prove (\ref{loss:pin-ball-more-a1}) and (\ref{loss:pin-ball-more-a2}),
we first consider the case  $t^*_{\mathrm{min}} = t^*_{\mathrm{max}}$. Then (\ref{loss:pin-ball-more-a1-h1}) and (\ref{loss:pin-ball-more-a1-h2})
yield $\CQ L t^*_{\mathrm{min}} = \CQ L t^*_{\mathrm{max}} \leq \CQ L t$, $t\in \mathbb{R}$. This implies
$\CQ L t^*_{\mathrm{min}} =  \CQ L t^*_{\mathrm{max}}  = \CQB L$,
and hence we conclude that (\ref{loss:pin-ball-more-a1-h1}) and (\ref{loss:pin-ball-more-a1-h2}) are equivalent to
(\ref{loss:pin-ball-more-a1}) and (\ref{loss:pin-ball-more-a2}), respectively.
Moreover, in the case $t^*_{\mathrm{min}} \neq  t^*_{\mathrm{max}}$, we have $\mathrm{Q}((t^*_{\mathrm{min}}(\mathrm{Q}), t^*_{\mathrm{max}}(\mathrm{Q}))) = 0$, which in turn
implies $\mathrm{Q}((-\infty,t^*_{\mathrm{min}}]) = \tau$ and $\mathrm{Q}([t^*_{\mathrm{max}},\infty)) = 1-\tau$.
For $t\in (t^*_{\mathrm{min}},t^*_{\mathrm{max}}]$, we consequently find
\begin{eqnarray}\label{h1}
 \CQ L t
&=&
(\tau-1)\int_{y<t}(y-t) \,  \mathrm{d}   \mathrm{Q}(y) + \tau \int_{y\geq t}(y-t)  \, \mathrm{d}  \mathrm{Q}(y)\nonumber
\\[-8pt]\\[-8pt]
& = &
(\tau-1)\int_{y<t^*_{\mathrm{max}}}y  \, \mathrm{d }  \mathrm{Q}(y) + \tau \int_{y\geq t^*_{\mathrm{max}}}y  \, \mathrm{d}  \mathrm{Q}(y)  ,\nonumber
\end{eqnarray}
where  we used $\mathrm{Q}((-\infty,t)) = \mathrm{Q}((-\infty,t^*_{\mathrm{min}}]) = \tau$ and  $\mathrm{Q}([t,\infty)) = \mathrm{Q}([t^*_{\mathrm{max}},\infty)) = 1-\tau$.
Since the right-hand side of (\ref{h1}) is independent of $t$, we thus conclude $\CQ L t = \CQ L t^*_{\mathrm{max}}$ for all
$t\in (t^*_{\mathrm{min}},t^*_{\mathrm{max}}]$. Analogously, we find $\CQ L t = \CQ L t^*_{\mathrm{min}}$ for all
$t\in [t^*_{\mathrm{min}},t^*_{\mathrm{max}})$, and hence we can, again, conclude $\CQ L t^*_{\mathrm{min}} = \CQ L t^*_{\mathrm{max}} \leq \CQ L t$ for all $t\in \mathbb{R}$.
As in the case $t^*_{\mathrm{min}} = t^*_{\mathrm{max}}$, the latter implies that (\ref{loss:pin-ball-more-a1-h1}) and (\ref{loss:pin-ball-more-a1-h2}) are equivalent to
(\ref{loss:pin-ball-more-a1}) and (\ref{loss:pin-ball-more-a2}), respectively.

For the proof of $\FQ L {0^+} = [t^*_{\mathrm{min}},t^*_{\mathrm{max}}]$, we first note that the previous discussion has already
shown $\FQ L {0^+} \supset [t^*_{\mathrm{min}},t^*_{\mathrm{max}}]$. Let us assume that $\FQ L {0^+} \not\subset [t^*_{\mathrm{min}},t^*_{\mathrm{max}}]$.
By a symmetry argument, we then may assume without loss of generality that there exists a $t\in \FQ L {0^+}$ with $t>t^*_{\mathrm{max}}$.
From (\ref{loss:pin-ball-more-a1}) we then conclude that $q_+ = 0$ and
$\mathrm{Q}((t^*_{\mathrm{max}}, t))=0$. Now, $q_+ = 0$ together with (\ref{loss:pin-ball-more-h1}) shows
$\mathrm{Q} ( (-\infty,t^*_{\mathrm{max}}] )  = \tau$, which in turn implies $\mathrm{Q} ( (-\infty,t] )  \geq \tau$.
Moreover, $\mathrm{Q}((t^*_{\mathrm{max}}, t))=0$ yields
\[
 \mathrm{Q}([t,\infty))
=
\mathrm{Q}([t^*_{\mathrm{max}},\infty)) - \mathrm{Q}(\{t^*_{\mathrm{max}}\})
=
1- \mathrm{Q} ( (-\infty,t^*_{\mathrm{max}}] )
= 1-\tau
  .
\]
 In other words, $t$ is a $\tau$-quantile, which contradicts  $t>t^*_{\mathrm{max}}$.
\end{pf}

For the proof of Theorem \ref{main1} we further need the
\textit{self-calibration loss} of  $L$ that is defined by
\begin{equation}\label{loss:self-cal-loss-def}
\breve{L}(\mathrm{Q},t) :=    \operatorname{dist}(t,\FQ L {0^+})  , \qquad t\in \mathbb{R}  ,
\end{equation}
where $\mathrm{Q}$ is a distribution   with $\operatorname{supp} \mathrm{Q}\subset [-1,1]$.
Let us  define  the \textit{self-calibration function}   by
\[
\dmaxexu {\varepsilon} {\mathrm{Q}}  {\breve{L},L} 
:=
   \inf_{t\in \mathbb{R}:   \breve{L}(\mathrm{Q},t) \ge \varepsilon} \CQ L t - \CQB L  ,  \qquad \varepsilon\geq 0.
\]
Note that if, for $t\in \mathbb{R}$, we write $\varepsilon:= \operatorname{dist}(t,\FQ L {0^+})$, then we have
 $\breve{L}(\mathrm{Q},t) \ge \varepsilon$, and hence the definition of the
self-calibration function yields
%
\begin{equation}\label{m1}
\dmaxexub {\operatorname{dist}(t,\FQ L {0^+})} {\mathrm{Q}}  {\breve{L},L} \leq \CQ L t - \CQB L  ,  \qquad t\in \mathbb{R}.
\end{equation}
In other words, the self-calibration function measures how well an  $\varepsilon$-approximate $L$-risk minimizer $t$ approximates
the set 
of exact $L$-risk
minimizers.

Our next goal is to estimate the self-calibration function for the pinball loss. To this end we need the following simple technical lemma.

\begin{lemma}\label{lower-pol}
For $\alpha\in [0,2]$ and $q\in [1,\infty)$ consider the function $\delta:[0,2]\to [0,\infty)$ defined by
\[
\delta(\varepsilon):=
\cases{
\varepsilon^q, &\quad \mbox{if }$  \varepsilon\in [0,\alpha]$,
\cr
q\alpha^{q-1} \varepsilon - \alpha^q(q-1), &\quad \mbox{if }$  \varepsilon\in [\alpha,2]$.
}
\]
Then, for all $\varepsilon\in [0,2]$, we have
\[
\delta(\varepsilon) \geq \biggl( \frac \alpha 2\biggr)^{q-1} \varepsilon^q  .
\]
\end{lemma}

\begin{pf}
Since $\alpha\leq 2$ and $q\geq 1$ we easily see by the definition of $\delta$ that the assertion is
true for $\varepsilon\in [0,\alpha]$. Now consider the function $h\dvtx [\alpha,2]\to \mathbb{R}$ defined by
\[
h(\varepsilon) := q\alpha^{q-1} \varepsilon - \alpha^q(q-1) - \biggl( \frac \alpha 2\biggr)^{q-1} \varepsilon^q  ,  \qquad \varepsilon\in [\alpha,2].
\]
It suffices
 to show that $h(\varepsilon)\geq 0$ for all $\varepsilon\in [\alpha,2]$. To show the latter
we first check that
\[
h'(\varepsilon) = q\alpha^{q-1}  - q\biggl( \frac \alpha 2\biggr)^{q-1} \varepsilon^{q-1}   ,  \qquad \varepsilon\in [\alpha,2]
\]
and hence we have $h'(\varepsilon)\geq 0$ for all $\varepsilon\in [\alpha,2]$.
Now we obtain the assertion from
this, $\alpha\in[0,2]$ and
\[
h(\alpha) = \alpha^q - \biggl( \frac \alpha 2\biggr)^{q-1} \alpha^q = \alpha^q \biggl( 1 -\biggl( \frac \alpha 2\biggr)^{q-1}\biggr)\geq 0  .
\]
\upqed\end{pf}

\begin{lemma}\label{self-cal-pinball-lower}
Let $L$ be the $\tau$-pinball loss and
$\mathrm{Q}$ be a distribution on $\mathbb{R}$ with $\operatorname{supp} \mathrm{Q}\subset [-1,1]$
that has a $\tau$-quantile  of type $q\in [1,\infty)$. Moreover,
let
$\alpha_\mathrm{Q}\in (0,2]$ and  $b_\mathrm{Q}>0$ denote the corresponding constants.
Then, for all $\varepsilon\in [0,2]$, we have
\[
\dmaxexu {\varepsilon} {\mathrm{Q}}  {\breve{L},L} \geq q^{-1}b_\mathrm{Q} \biggl(\frac {\alpha_\mathrm{Q}}2 \biggr)^{q-1} \varepsilon^q =
 q^{-1}2^{1-q} \gamma_\mathrm{Q} \varepsilon^q          .
\]
\end{lemma}

\begin{pf}
Since $L$ is convex, the map
$t\mapsto \CQ L  t - \CQB L$ is convex, and thus it is decreasing
on $(-\infty, t^*_{\mathrm{min}}]$ and
increasing on $[t^*_{\mathrm{max}},\infty)$.
Using $\FQ L {0^{+}} = [t^*_{\mathrm{min}},t^*_{\mathrm{max}}]$, we thus find
\[
\FQ {\breve{L}} {\varepsilon} := \{t\in \mathbb{R}\dvtx \breve{L}(\mathrm{Q},t)<\varepsilon\}
= (t^*_{\mathrm{min}}-\varepsilon,t^*_{\mathrm{max}}+\varepsilon)
\]
for all $\varepsilon>0$. Since
this gives $\dmaxexu {\varepsilon} {\mathrm{Q}}  {\breve{L},L}\nonumber
=
\inf_{t\notin \FQ {\breve{L}} {\varepsilon}} \CQ L t-\CQB L$, we obtain
\begin{equation}\label{loss:self-cal-conv-ineq}
\dmaxexu {\varepsilon} {\mathrm{Q}}  {\breve{L},L}
=
\min\{  \CQ L {t^*_{\mathrm{min}}-\varepsilon}, \CQ L {t^*_{\mathrm{max}}+\varepsilon}      \} - \CQB L   .
\end{equation}
Let us first consider the case $q\in (1,\infty)$. For $\varepsilon\in [0,\alpha_\mathrm{Q}]$,
(\ref{loss:pin-ball-more-a1}) and (\ref{q-type-2}) then yield
\[
\CQ L {t^*_{\mathrm{max}}+\varepsilon} - \CQB L
=
 \varepsilon q_++\int_0^\varepsilon \mathrm{Q}\bigl( (t^*_{\mathrm{max}},t^*_{\mathrm{max}}+s)\bigr) \,  \mathrm{d}s
\geq
b_\mathrm{Q}  \int_0^\varepsilon  s^{q-1}  \, \mathrm{d}s
 =
q^{-1} b_\mathrm{Q} \varepsilon^q  ,
\]
and, for $\varepsilon\in [\alpha_\mathrm{Q},2]$, (\ref{loss:pin-ball-more-a1}) and (\ref{q-type-2})  yield
\[
\CQ L {t^*_{\mathrm{max}}+\varepsilon} - \CQB L
 \geq
b_\mathrm{Q}  \int_0^{\alpha_\mathrm{Q}}  s^{q-1} \,  \mathrm{d}s + b_\mathrm{Q}  \int_{\alpha_\mathrm{Q}}^\varepsilon  \alpha_\mathrm{Q}^{q-1}  \, \mathrm{d}s
 =
q^{-1} b_\mathrm{Q} \bigl( q\alpha_\mathrm{Q}^{q-1} \varepsilon - \alpha_\mathrm{Q}^q(q-1) \bigr)   .
\]
For $\varepsilon\in [0,2]$, we have thus shown $\CQ L {t^*_{\mathrm{max}}+\varepsilon} - \CQB L  \geq q^{-1} b_\mathrm{Q} \delta(\varepsilon)$, where $\delta$ is the function defined in Lemma
\ref{lower-pol} for $\alpha:= \alpha_\mathrm{Q}$.

Furthermore, in the case $q=1$ and $t^*_{\mathrm{min}} \neq t^*_{\mathrm{max}}$, Proposition \ref{loss:pin-ball-more}
shows $q_+ = \mathrm{Q}(\{t^*_{\mathrm{max}}\})$, and hence (\ref{loss:pin-ball-more-a1})
yields $\CQ L {t^*_{\mathrm{max}}+\varepsilon} - \CQB L \geq \varepsilon q_+ \geq b_\mathrm{Q} \varepsilon$ for all
$\varepsilon\in [0,2] = [0,\alpha_\mathrm{Q}]$ by the definition of $b_\mathrm{Q}$ and $\alpha_\mathrm{Q}$.
In the case $q=1$ and $t^*_{\mathrm{min}} = t^*_{\mathrm{max}}$, (\ref{loss:pin-ball-more-h1}) yields
$q_+ = \mathrm{Q}((-\infty,t^*]) - \tau \geq b_\mathrm{Q}$ by the definition of $b_\mathrm{Q}$, and hence (\ref{loss:pin-ball-more-a1})
again gives $\CQ L {t^*_{\mathrm{max}}+\varepsilon} - \CQB L \geq  b_\mathrm{Q} \varepsilon$ for all $\varepsilon\in [0,2]$.
Finally, using (\ref{loss:pin-ball-more-a2}) instead of (\ref{loss:pin-ball-more-a1}),
we can analogously show
$\CQ L {t^*_{\mathrm{min}}-\varepsilon} - \CQB L  \geq q^{-1} b_\mathrm{Q} \delta(\varepsilon)$ for all $\varepsilon\in  [0,2]$ and $q\geq 1$.
By
(\ref{loss:self-cal-conv-ineq}) we thus conclude
that
\[
\dmaxexu {\varepsilon} {\mathrm{Q}}  {\breve{L},L} \geq q^{-1} b_\mathrm{Q} \delta(\varepsilon)
\]
for all $\varepsilon\in [0,2]$. Now the assertion follows from Lemma \ref{lower-pol}.
\end{pf}

\begin{pf*}{Proof of Theorem \ref{main1}}
For fixed $x\in X$ we write $\varepsilon := \operatorname{dist}(f(x),\FPo L {0^+})$.
By Lemma~\ref{self-cal-pinball-lower} and (\ref{m1}) we obtain, for $\mathrm{P}_X$-almost all $x\in X$,
\begin{eqnarray*}
\bigl|\operatorname{dist}\bigl(f(x),\FPo L {0^+}\bigr)\bigr|^q
&\leq &
q 2^{q-1} \gamma^{-1}(x) \dmaxexub {\varepsilon} {\mathrm{P}(\cdot|x)}  {\breve{L},L}\\
& \leq & q 2^{q-1}\gamma^{-1}(x) \bigl(  \CPo L (f(x)) - \CPoB L\bigr)  .
\end{eqnarray*}
By taking the $\frac p {p+1}$th power on both sides, integrating and finally
applying H\"{o}lder's inequality, we then obtain the assertion.
\end{pf*}

\begin{pf*}{Proof of Theorem \ref{main2}}
Let $f\dvtx X\to [-1,1]$ be a  function.
Since $F^*_{\tau,\mathrm{P}}(x)$ is closed, there then exists a $\mathrm{P}_X$-almost surely uniquely determined
function $f^*_{\tau,\mathrm{P}}\dvtx X\to [-1,1]$
that satisfies both
\begin{eqnarray*}
f^*_{\tau,\mathrm{P}}(x) & \in & F^*_{\tau,\mathrm{P}}(x),
\\
|f(x) - f^*_{\tau,\mathrm{P}}(x)| &=& \operatorname{dist}( f(x), F^*_{\tau,\mathrm{P}}(x))
\end{eqnarray*}
for $\mathrm{P}_X$-almost all $x\in X$.
Let us write $r:= \frac {pq} {p+1}$.
We first consider the case $r\leq 2$, that is, $\frac 2 q \leq \frac p {p+1}$.
Using the Lipschitz continuity of the pinball loss $L$ and Theorem \ref{main1} we then obtain
\begin{eqnarray*}
\mathbb{E}_\mathrm{P} (  L\circ f - L\circ \fpb     )^2
&\leq &
\mathbb{E}_{\mathrm{P}_X} | f -\fpb  |^2
\\
&\leq &
\snorm{f -\fpb}_\infty^{2-r}  \mathbb{E}_{\mathrm{P}_X} | f -\fpb  |^r
\\
&\leq &
2^{2-r/q} q^{r/q}
\snorm{\gamma^{-1}}_{L_p(\mathrm{P}_X)}^{r/q}  \bigl(\RP L f - \RPB L \bigr)^{r/q}  .
\end{eqnarray*}
Since $\frac{r}{q}  =\frac p {p+1} = \vartheta$, we thus obtain the assertion in this case.
Let us now consider the case $r>2$. The
Lipschitz continuity of $L$ and Theorem \ref{main1} yield
\begin{eqnarray*}
\mathbb{E}_\mathrm{P} (  L\circ f - L\circ \fpb     )^2
&\leq &
\bigl(\mathbb{E}_\mathrm{P} (  L\circ f - L\circ \fpb     )^r\bigr)^{2/r}
\\
&\leq &
(\mathbb{E}_{\mathrm{P}_X} | f -\fpb  |^r )^{2/r}
\\
&\leq &
\bigl(2^{1-1/q} q^{1/q}  \snorm{\gamma^{-1}}_{L_p(\mathrm{P}_X)}^{1/q}  \bigl(\RP {L} f - \RPB {L}  \bigr)^{1/q} \bigr)^{2}
\\
& =  &
2^{2-2/q} q^{2/q}  \snorm{\gamma^{-1}}_{L_p(\mathrm{P}_X)}^{2/q}  \bigl(\RP {L} f - \RPB {L}  \bigr)^{2/q}   .
\end{eqnarray*}
Since for $r>2$ we have $\vartheta = 2/q$, we again obtain the assertion.
\end{pf*}

\begin{pf*}{Proof of Theorem \ref{oracle-general}}
As shown in \cite{Ste08a}, Lemma 2.2,    (\ref{eigenvalues-assump}) is equivalent to
the entropy assumption (\ref{EN}), which in turn implies (see \cite{Ste08a}, Theorem 2.1, and
\cite{SteinwartChristmann2008a}, Corollary 7.31)
\begin{equation}\label{EEN}
\mathbb{E}_{D_X\sim \mathrm{P}_X^n} e_i\bigl(\operatorname{id}\dvtx H\to \Lx 2 {\mathrm{D}_X}\bigr) \leq  c \sqrt a  i^{-1/(2\varrho)} ,  \qquad i\geq 1  ,
\end{equation}
where $\mathrm{D}_X$ denotes the empirical measure with respect to $D_X=(x_1,\dots,x_n)$ and $c\geq 1$ is
a constant only depending on $\varrho$. Now the assertion follows from \cite{SteinwartChristmann2008a}, Theorem 7.23,
by considering the  function $f_0\in H$ that achieves
$\lambda \snorm {f_0}_H^2 + \RP L {f_0} - \RPB L = A(\lambda)$.
\end{pf*}

\printhistory

\end{document}